\documentclass[12pt,a4paper,twoside]{article}
\usepackage{amsthm,amsmath,amssymb}
\newtheorem{teor}{Theorem}[section]
\newtheorem{defin}[teor]{Definition}
\newtheorem{lemm}[teor]{Lemma}
\newtheorem{osse}[teor]{Remark}
\newtheorem{prop}[teor]{Proposition}
\newtheorem{defi}[teor]{Definition}
\newtheorem{coro}[teor]{Corollary}
\newtheorem{prob}[teor]{Problem}

\newcommand{\bele}{\begin{lemm}\begin{sl}}
\newcommand{\enle}{\end{sl}\end{lemm}}
\newcommand{\bedef}{\begin{defi}\begin{sl}}
\newcommand{\eddef}{\end{sl}\end{defi}}
\newcommand{\bete}{\begin{teor}\begin{sl}}
\newcommand{\ente}{\end{sl}\end{teor}}
\newcommand{\beos}{\begin{osse}\begin{rm}}
\newcommand{\eddos}{\end{rm}\end{osse}}
\newcommand{\bepr}{\begin{prop}\begin{sl}}
\newcommand{\empr}{\end{sl}\end{prop}}
\newcommand{\bepro}{\begin{prob}\begin{rm}}
\newcommand{\empro}{\end{rm}\end{prob}}
\newcommand{\bede}{\begin{defin}\begin{sl}}
\newcommand{\edde}{\end{sl}\end{defin}}
\newcommand{\beco}{\begin{coro}\begin{sl}}
\newcommand{\enco}{\end{sl}\end{coro}}

\newcommand{\quext}{\quad\text}

\newcommand{\de}{\partial}

\newcommand{\RR}{\mathbb{R}}
\newcommand{\NN}{\mathbb{N}}

\newcommand{\beeq}[1]{\begin{equation}\label{#1}}
\newcommand{\eddeq}{\end{equation}}

\newcommand{\beeqa}[1]{\begin{eqnarray}\label{#1}}
\newcommand{\eddeqa}{\end{eqnarray}}

\newcommand{\beal}[1]{\begin{align}\label{#1}}
\newcommand{\eddal}{\end{align}}

\newcommand{\bespl}[1]{\begin{split}\label{#1}}
\newcommand{\edspl}{\end{split}}

\newcommand{\bega}[1]{\begin{gather}\label{#1}}
\newcommand{\edga}{\end{gather}}

\newcommand{\beeqax}{\begin{eqnarray*}}
\newcommand{\eddeqax}{\end{eqnarray*}}

\def\qed{\ifmmode 
  \else \leavevmode\unskip\penalty9999 \hbox{}\nobreak\hfill
  \fi
  \quad\hbox{\hskip.5em\vrule width.4em height.6em depth.05em\hskip.1em}}
\def\endproofsym{\qed}
\newcommand{\dimbox}{\hbox{\hskip.5em\vrule width.4em height.6em depth.05em\hskip.1em}}

\def\endnobox{\def\endproofsym{}\end{proof}\def\endproofsym{\qed}}

\newcommand{\no}{\nonumber}

\newcommand{\beeqao}{\begin{eqnarray}\no}
\newcommand{\bealo}{\begin{align}\no}
\newcommand{\besplo}{\begin{split}\no}
\newcommand{\begao}{\begin{gather}\no}

\newcommand\cpt{\hskip1.15truept\cdot
               \hskip1.15truept} 

\newcommand{\nor}[2]{\|#1\|_{#2}}
\newcommand{\nord}[2]{\|#1\|^2_{#2}}

\newcommand{\unmezzo}{\frac{1}{2}}

\newcommand{\dua}[3]{{}_{#2}\langle{#1}\rangle_{#3}}

\newcommand{\+}{\hspace{1pt}}
\newcommand{\perogni}{\forall\,}
\newcommand{\esiste}{\exists\,}
\newcommand{\itt}{\int_0^t}

\newcommand{\io}{\int_\Omega}

\newcommand{\iga}{\int_\Gamma}

\newcommand{\epsi}{\varepsilon}

\newcommand{\lhs}{left hand side}
\newcommand{\rhs}{right hand side}

\DeclareMathOperator{\deriv}{d}

\DeclareMathOperator{\dom}{dom}

\DeclareMathOperator{\Id}{Id}

\DeclareMathOperator{\sign}{sign}

\DeclareMathOperator{\loc}{loc}

\let\TeXchi\chi
\def\chi{{\setbox0 \hbox{\mathsurround0pt
$\TeXchi$}\hbox{\raise\dp0 \copy0 }}}

\newcommand{\teta}{\vartheta}

\newcommand{\tetan}{\teta_n}
\newcommand{\chin}{\chi_n}
\newcommand{\chint}{\chi_{n,t}}

\newcommand{\calX}{{\cal X}}

\newcommand{\calE}{{\cal E}}
\newcommand{\calF}{{\cal F}}

\newcommand{\calY}{{\cal Y}}
\newcommand{\calZ}{{\cal Z}}

\newcommand{\calJ}{{\cal J}}

\newcommand{\calV}{{\cal V}}

\newcommand{\barO}{\overline{\Omega}}

\newcommand{\chii}{\chi_\infty}

\newcommand{\dit}{\deriv\!t}
\newcommand{\dis}{\deriv\!s}
\newcommand{\dix}{\deriv\!x}

\newcommand{\ddt}{\frac{\deriv\!{}}{\dit}}

\numberwithin{equation}{section}

\begin{document}
\thispagestyle{empty}
\begin{center}

{\LARGE Long time convergence for a class\\[0.2truecm]
of variational phase field models\footnotemark[1]}
\footnotetext[1]{This paper was initiated during a visit of
the first author to the Universit\'e Paris-Sud XI,
Laboratoire d'Analyse Num\'erique et EDP,
whose kind hospitality is gratefully acknowledged.
The work also benefited from a financial support from
the MIUR-COFIN 2004 research program on ``Mathematical
modelling and analysis of free boundary problems''.}

\font\autore=cmcsc10                    
\font\address=cmsl10 at 10.95 truept    
\font\email=cmtt10 at 10.95truept       
\font\auts=cmcsc10 at 10.95 truept      
\font\sls=cmsl10 at 10.95 truept        
\font\bfs=cmbx10 at 10.95 truept        
\font\smallfont=cmr10 at 10.95 truept   
\def\Email{\smallfont E-mail: \email}

\vskip.75truecm

\autore
Pierluigi Colli$^{(1)}$

\Email pierluigi.colli@unipv.it
\vskip0.5truecm

\autore
Danielle Hilhorst$^{(2)}$

\Email danielle.hilhorst@math.u-psud.fr
\vskip0.5truecm

\autore
 Fran\c coise Issard-Roch$^{(2)}$

\Email francoise.issard-roch@math.u-psud.fr
\vskip0.5truecm

\autore
Giulio Schimperna$^{(1)}$

\Email giusch04@unipv.it

\vskip0.75truecm

\address
$^{(1)}$ Dipartimento di Matematica ``F. Casorati''

Universit\`a di Pavia, Via Ferrata 1, I-27100 Pavia, Italy

\vskip0.5truecm

$^{(2)}$ CNRS and Laboratoire de Math\'ematiques

Universit\'e Paris-Sud 11, Bat.~425, F-91405 Orsay, France

\end{center}

\noindent
\centerline{\it Dedicated to Professor Masayasu Mimura on the occasion of his 65th birthday}

\vskip.75truecm\noindent
%
%
%
%
%
{\small\bf Abstract.} {\small\rm In this paper we analyze
a class of phase field models for the dynamics of phase transitions
which extend the well-known Caginalp and Penrose-Fife models.
Existence and uniqueness of the solution
to the related initial boundary value problem are shown.
Further regularity of the solution is deduced by
exploiting the so-called regularizing effect.
Then, the large time behavior of such a solution is
studied and several convergence properties of the trajectory
as time tends to infinity are discussed.}

\vspace{.2cm}


\vspace{2mm}

\noindent
{\small\bf AMS (MOS) subject clas\-si\-fi\-ca\-tion:}
{\small\rm 35B40, 35K45, 80A22.}

\vspace{2mm}



\pagestyle{myheadings}
\newcommand\testopari{\sc Colli -- Hilhorst -- Issard-Roch -- Schimperna}
\newcommand\testodispari{\sc long time convergence for phase field models}
\markboth{\testodispari}{\testopari}


\section{Introduction}
\label{intro}

The present note is devoted to the analysis of the
regularity and long-time behavior properties of the
following class of PDE's systems modeling phase
change phenomena
\beal{calorein}
   & \epsi\+\teta_t+\lambda(\chi)_t-\Delta j'(\teta)=f,\\
  \label{phasein}
   & \delta\+\chi_t-\Delta\chi+W'(\chi)=
    \lambda'(\chi)\+j'(\teta).
\end{align}
Here, the unknowns are the relative temperature $\teta$
(i.e., some critical freezing or melting
temperature $\tau_c$ has been normalized to 0)
and the order parameter, or phase field, $\chi$;
both are functions of the spatial variable $x$
(ranging in a bounded, connected, and sufficiently smooth
domain $\Omega\subset\RR^d$, $1\le d\le 3$) and of the
time $t\in[0,\infty)$ (let us use the shorter notation $\infty$ in place of
$+\infty$). In the above system, $\epsi,\,\delta>0$ are relaxation parameters;
the function $W:\/\dom W\to \RR$ represents a configuration potential in $\chi$
and is assumed to be convex in its principal part; $\lambda(\cdot)$
is a possibly non linear function with $\epsi \teta + \lambda(\chi)$
yielding the {\sl internal energy}\/ of the system;
finally, $f$ stands for the heat source.
However, the main novelty in \eqref{calorein}--\eqref{phasein} is
given by the presence of the nonlinear, but convex
function $j$. This nonlinearity can be justified
both on the physical and on the analytical side.
Actually, several well-known models can
be included in this general framework.
For example, $j(r)=r^2/2$, corresponds to
the Caginalp system \cite{Cag},
while $j(r)=-\log(r+\tau_c)+r/\tau_c$
gives the so-called Penrose-Fife model \cite{PF,PF2}
(the complication of the expression, which is not
the usual one, is compensated by the nice property $j'(0)=0 $).
Furthermore, also intermediate choices for
$j$ correspond to meaningful cases:
for instance, a combination of the previous
two expressions provides
a variant of the Penrose-Fife model
with {\sl special}\/ heat flux law introduced in
\cite{CL,CLS}.

Various systems concerned with the abovementioned phase field models,
possibly including non-smooth potentials and also applying to
martensitic phase transformations, have been intensively investigated in the last years.
Among a number of recent contributions, let us quote
\cite{IKK} and \cite{CGRS} which treat the case of Neumann boundary
conditions for the temperature; \cite{KK} devising a general (convex) framework for the study of
Penrose-Fife systems; \cite{CP} involved with the analysis of the quasistationary (i.e.,
$\delta =0$ in \eqref{phasein}) Penrose-Fife model;
\cite{IK} that deals with the long-time behaviour and the study of inertial sets;
\cite{AF} and \cite{AFIR} studying the long-time convergence of the Caginalp model with and
without memory effects, and \cite{SZ, FSc} addressing the same questions for
the Penrose-Fife model;  \cite{RS, RS4} showing the existence of a uniform attractor
for Penrose-Fife systems; \cite{CGI} and \cite{RS6} in which a
hyperbolic dynamics for $\chi$, characterized by an extra inertial
term $\rho  \chi_{tt}$ in the \lhs\ of equation \eqref{phasein}, is considered.

In all this framework, the occurrence of a {\sl general}\/ convex function $j$
seems to be an interesting isssue, worth to be deepened. Recently,
well-posedness of an initial and boundary value problem
for \eqref{calorein}--\eqref{phasein}
has been shown in \cite{Ro} for completely
arbitrary (convex) $j$. More precisely, in \cite{Ro}
the homogeneous Neumann problem related to
\eqref{calorein}--\eqref{phasein} is studied
for a {\sl linear}\/ function $\lambda$.
We remark that, while the Neumann boundary conditions
for the phase field $\chi$ appear to be the most
natural choice for phase field models,
the case of no-flux conditions for $j'(\teta)$ was
motivated in \cite{Ro} by the purpose of studying
some singular limits of the system. In particolar, in \cite{Ro}
it is shown that \eqref{calorein}--\eqref{phasein}
gives rise to the Cahn-Hilliard equation
in the {\sl viscous}\/ form
if $\epsi$ is sent to 0 and in the
standard form if both $\epsi,\delta$ go to 0.
On the other hand, the choice of the no-flux
conditions for $j'(\teta)$ gives rise to some additional
difficulties in the analysis that forced the author
of \cite{Ro} to restrict the range of the admissible
potentials $W$ in order to get some a priori control of the
spatial average of the unknowns
(see \cite[assumption (28)]{Ro}, see also \cite{CGRS}).

In this paper, we first extend the well-posedness result of \cite{Ro},
by adapting it to our slightly different setting. Then,
we investigate here some further
properties of the solution to suitable initial-boundary
value problems related to \eqref{calorein}--\eqref{phasein}.
Namely, we shall concentrate our attention on the long-time behavior of
the system from the point of view of $\omega$-limits
of solution trajectories. Since we are not interested
in singular limits, {\sl we shall take $\epsi=\delta=1$
in the sequel}.

In our approach, the basic observation is that
the system, at least in the case of no external source,
admits the Liapounov functional
\beeq{energyin}
  \calE(\teta,\chi):=\io\Big(
   \frac12|\nabla\chi|^2
   +W(\chi)
   +j(\teta)\Big),
\end{equation}
which is obtained testing \eqref{calorein} by $j'(\teta)$ and
\eqref{phasein} by $\chi_t$, then taking the sum.
Hence, if we introduce the new variable $e:=\teta+\lambda(\chi)$
(internal energy), \eqref{calorein}--\eqref{phasein}
can be seen (at least in the case of no-flux conditions and no
external source) as a {\sl gradient flow}\/ problem
\beeq{gfin}
  e_t=-\de_{*,e}\calE, \qquad
  \chi_t=-\de_{H,\TeXchi}\calE.
\end{equation}
Here, the symbol $\de$ denotes here (sub)differentiation, and,
more precisely, $\de_{H,\TeXchi}$ stands for the subdifferential
w.r.t.~$\chi$ in $H:=L^2(\Omega)$ while $\de_{*,e}$
indicates differentiation in $e$ in the space $H^1(\Omega)^*$
(which gives rise to the Laplacian in \eqref{calorein}).

In view of this variational structure, it is reasonable
to expect good asymptotic properties of the solution
as $t$ goes to infinity. This has already been noticed in some
particular cases. For instance, for the Caginalp model
corresponding to $j(r)=r^2/2$, the long time analysis
has been performed in \cite{AFIR,GPS1} for various types of
boundary conditions and assumptions on $\lambda,\, W$.
Instead, the Penrose-Fife case ($j(r)=-\log r$) has been studied
in \cite{FSc}, referring only to
the case of (nonhomogeneous) Dirichlet
boundary conditions for the temperature. Concerning
\cite{FSc}, we point out that the problem
addressed there presents a number of
mathematical difficulties due the character
of $j(\cdot)=-\log(\cdot)$, which is both
singular at $0$ and unbounded from below.
Actually, the precompactness of trajectories, leading to
the existence of a nonempty $\omega$-limit set,
has been proved in \cite{FSc} by strongly
relying on {\sl ad hoc}\/ techniques to overcome,
in particular, the noncoercive character of $j$ at $+\infty$.

Based on these considerations, it seems very difficult
to perform a long time analysis of \eqref{calorein}--\eqref{phasein}
by taking $j$ completely general as in \cite{Ro}. Rather, we
find there is room for a bit of compromise and ask
at least that $j$ has some coercivity: see \eqref{hpj} below.
On the other hand, $j$ can exhibit a {\sl singular}\/
character too. Still for the purpose of coercivity,
we also need boundary conditions for $\teta$ which give some
(uniform in time) control on its space average
(differently from the no-flux case of \cite{Ro}).
Then, we shall consider two sets of boundary conditions:
\beal{DBC}
  & i)\quad \teta=\teta_\infty,~~\hbox{with }\, j'(\teta_\infty)=0,~~
   \hbox{and }\, \partial_n\chi=0,
\end{align}
which we refer to as Dirichlet boundary conditions, and
\beal{RBCin}
   &i i)\quad -\partial_n j'(\teta)=\eta (j'(\teta)-j'(\teta_\Gamma))
    ~~\hbox{and }\,\partial_n\chi=0,
\end{align}
which we refer to as Robin boundary conditions, where
$\teta_\Gamma$ represents the extremal boundary temperature and
$\eta$ is a positive constant. We point out that
the latter case subsumes the presence of some (boundary)
source term. We will see that this does not destroy the
variational character of the system, nor does this the presence of a
nonzero volumetric heat source $f$ in \eqref{calorein},
provided $f$ and the boundary datum are globally $L^2$ in time.
In our analysis, we will be able to
consider nonlinear latent heat functions $\lambda$ and,
what is more important, rather general potentials $W$
with the only restriction that they should not
exhibit minima at the boundary of their domain
(like instead it could happen in the case of a
double-obstacle potential).
These could neither be considered in \cite{FSc},
essentially for technical
reasons, nor in \cite{Ro} due to the quoted difficulties
coming from the no-flux conditions for $j'(\teta)$.
We stress in particular that
{\sl singular}\/ potentials $W$ (i.e.~which are $+\infty$ outside
an interval $I\subset\RR$) are not completely easy
to address even in the simpler
case of the Caginalp model \cite{GPS1}.

In the framework of the long time analysis, we prove the following
results. First, we demonstrate that the $\omega$-limit of any solution
trajectory is not empty and consists only of steady-state solutions
$(\teta_\infty , \chi_\infty) $. More in detail, the component
$\teta_\infty$ coincides with the (constant) temperature on
the boundary (i.e., with the {\sl unique}\/ value
such that $j'(\teta_\infty)=0$, cf.~assumption
\eqref{hpj2} below); so that we also obtain the
limit $\teta (t) \to \teta_\infty$ holds
in a suitable sense as $t\nearrow \infty $
and not only for a subsequence.
Conversely, in general we only have that,
as $t\nearrow\infty$, 
$\chi (t)$ is precompact in a suitable
topology and any of its limit points is a steady
state solution of \eqref{phasein}. In fact,
\eqref{phasein} may have infinitely many
stationary states due to nonconvexity of $W$.
We also point out that a careful use
of parabolic regularization effects is
a key step in the proof of the convergence result.
Our second theorem, which follows the lines
of some recent work devoted to the large time behavior
of degenerate parabolic equations~\cite{FS}
and phase transition systems (see, among others, \cite{AFIR,FIRP2,FSc}),
gives a sufficient condition under which
the $\omega$-limit consists of only one point.
Namely, we can prove that, in the case when
$W$ is analytic in the subdomain where the
solution component $\chi$ lives, then the $\omega$-limit set
is a singleton and consequently also the entire trajectory of $\chi(t)$
converges to $\chi_\infty.$ Here, the main ingredient of the proof
is the so-called {\sl Simon-\L ojasiewicz}\/ inequality \cite{Lo, Si2}, that
was originally \cite{Lo} stated as a nontrivial (local)
growth estimate for analytic functions of several
complex variables in the neighbourhood of a critical point:
Simon contributed by extending this inequality to the infinite dimensional
setting, thus allowing to characterize  the large time behavior of evolution systems
with {\sl analytic}\/ nonlinearities. The last result we present is a convergence
result which also establishes a (say, polynomial) rate of
convergence for the $L^2(\Omega)$-norm of $\chi(t) - \chi_\infty$.
In proving such a result, we follow a method from \cite{CJ}
and argue partly as in \cite{GPS2}.

Here is the plan of the paper. In  Section~\ref{mainres} we introduce
the functional framework, state precise hypotheses on the
data and formulate our main results.
Section~\ref{secunif} is devoted to the proof of the well-posedness of
the system (Theorem~\ref{teoexist})
and of the basic uniform estimates (Theorem~\ref{teorego}).
Finally, Section~\ref{secomegalim} is concerned
with all the properties of the $\omega$-limit set (proof of
Theorem~\ref{teoomega1} and Theorem~\ref{teoomega2}).


\section{Main results}
\label{mainres}

Let $\Omega$ be a $C^{1,1}$, bounded, and connected domain
in $\RR^d$, $1\le d\le 3$, and let $\Gamma:=\de\Omega$.
Set $H:=L^2(\Omega)$, $V:=H^1(\Omega)$, both endowed
with their standard scalar products and norms.
The norms in $H$ and in $H^d$ will be equally indicated
by $|\cdot|$ and $(\cdot,\cdot)$ will denote the
corresponding scalar products. Let also
$V_0:=H^1_0(\Omega)$, endowed with the norm
$\|\cdot\|_{V_0}:=|\nabla\cdot|$.
The symbol $\|\cdot\|_X$ will stand
for the norm in the generic Banach space $X$
and $\dua{\cdot,\cdot}{X^*}{X}$ will denote
the duality between $X$ and its topological
dual space $X^*$. The space $H$ will be identified
with its dual, so that we have the chains
of continuous embeddings $V\subset H\subset V^*$
and $V_0\subset H\subset V_0^*$.
We introduce the elliptic operator
\beeq{defiA}
  A:V\to V^*,\qquad
   \dua{Av,z}{V^*}V:=\io\nabla v\cdot\nabla z,
   \quext{for }\,v,z\in V.
\end{equation}
We also let $I,J$ be {\sl open}\/ intervals of $\RR$,
with $0\in I$, and let $I_0$ be an {\sl
open and bounded}\/ interval containing
$0$ and whose closure is contained in $I$.
Then, our basic hypotheses on the data are
\beal{hplambda}\tag{hp{$\lambda$}}
  & \lambda\in C^{1,1}_{\loc}(\RR),\quad
   \esiste \Lambda>0:~|\lambda''(r)|\le \Lambda~\,
   \text{for a.e.~}\,r\in\RR,\\
 \label{hpW}\tag{hp{$W$}\/1}
  & W\in C^{1,1}_{\loc}(I;[0,+\infty)),\quad
   \esiste \kappa>0:~W''(r)\ge -\kappa~\,
   \text{for a.e.~}\, r\in I,\\
 \label{hpW2}\tag{hp{$W$}\/2}
  & \esiste\mu>0:
   W'(r)/r\ge \mu~\,\text{for a.e.~}\,
    r\in I\setminus I_0,\\
 \label{hpj}\tag{hp{$j$}\/1}
  & j\in C^{1,1}_{\loc}(J;[0,+\infty)),\quad
   \esiste\sigma>0:~j''(r)\ge\sigma~\,\text{for a.e.~}\,r\in J,\\
 \label{hpj2}\tag{hp{$j$}\/2}
  & \esiste\teta_\infty\in J:~j'(\teta_\infty)=0.
\end{align}
Let us spend some words on the hypotheses on $W$ and $j$.
Formula \eqref{hpW} says that $W$ can be nonconvex, but
just up to a quadratic perturbation (actually,
$r\mapsto W(r)+\kappa r^2/2$ is convex).
Instead, $j$ is uniformly strictly convex by \eqref{hpj};
moreover, it has some coercivity property in the sense that,
by \eqref{hpj2}, it attains its minimum value
(which can be fixed at $0$, for simplicity)
at the point $\teta_\infty\in J$. The role
of \eqref{hpW2} will be outlined later on.

In the sequel, both $W$ and $j$ will be extended to the whole
real line by means of the following procedure.
First, we prolong $j$ (resp.~$W+\kappa\Id^2/2$) giving it
the value $+\infty$ outside $J$ (resp.~$I$); then, we take the
lower semicontinuous regularization, and,
finally, only from $W$ we subtract $\kappa\Id^2/2$.
This means, in the simpler (because convex) case of $j$
that, if $J$ is bounded and $j$ ``explodes''
(taking the value $+\infty$)
at its boundary, we simply extend it at $+\infty$ outside
$J$. If, instead, $J$ is bounded but $j$ does not explode
at least on one side of $J$, then we first close the graph
of $j$ and then extend it at $+\infty$. For $W$ the procedure
is slightly more complicated due to its possibly
non convex character. In any case, the extended
$j$ and $W+\kappa\Id^2/2$ are convex and lower
semicontinuous functions from $\RR$ to $[0,+\infty]$,
so that their $\RR$-subdifferentials
are {\sl maximal monotone graphs}\/ coinciding,
respectively on $J$ and on $I$,
with the ``original'' functions $j'$ and $W'+\kappa\Id$.
%
%

With all these conventions in mind,
our assumptions on the initial data are
\beal{hpteta0}\tag{hp{$\teta_0$}}
  & \teta_0\in H,\qquad
   j(\teta_0)\in L^1(\Omega),\\
 \label{hpchi0}\tag{hp{$\chi_0$}}
  & \chi_0\in V, \qquad
   W(\chi_0)\in L^1(\Omega).
\end{align}
\beos
 Recalling \eqref{energyin}, let us notice that
 \eqref{hpteta0}--\eqref{hpchi0} are equivalent
 to asking $\calE(\teta_0,\chi_0)<+\infty$,
 i.e., that the initial energy is finite. Indeed, the property
$\teta_0\in H $ follows from the quadratic growth of $j$ (cf.~\eqref{hpj}).
 Note also that, by \eqref{hpW} and \eqref{hpj}, $W$ and $j$ assume only non-negative values and
 consequently the functional $\calE$ is non-negative.
\eddos
\noindent%
Concerning the heat source, we assume in any case
\beeq{f1}\tag{hp{$f$}\/}
  f\in L^2(0,\infty;H),
\end{equation}
but this hypothesis will need some refinement in the sequel.

In case we work with the Robin boundary condition
\beal{RBC}
  & -\partial_n j'(\teta)=\eta (j'(\teta)-j'(\teta_\Gamma))~~
   \hbox{and }~\partial_n\chi=0,
\end{align}
where $\eta$ is some positive constant,
we also suppose that the external
boundary temperature $\teta_\Gamma$ satisfies
\beeq{tetaga}\tag{hp{$\teta_\Gamma$}}
  \teta_\Gamma:\Gamma\times(0,\infty)\to\RR~\,\text{measurable,}\qquad
   j'(\teta_\Gamma)\in L^2(0,\infty;L^2(\Gamma)).
\end{equation}
We remark that, by \eqref{hpj} (which ensures the
Lipschitz continuity of $(j')^{-1}$) and \eqref{hpj2},
\eqref{tetaga} entails that, for $n\nearrow\infty$
($n$ denoting here the time variable),
\beeq{bordazero}
  j'(\teta_\Gamma(n+\cdot))\to 0,\quad
   \teta_\Gamma(n+\cdot)\to \teta_\infty
   \quext{strongly in }\,L^2(0,1;L^2(\Gamma)).
\end{equation}
Next, we introduce the operator
\beeq{defiR}
  R:V\to V^*,\qquad
   \dua{Rv,z}{V^*}V=\io \nabla v\cdot\nabla z
    +\eta\iga v\+z.
\end{equation}
Of course, $R$ turns out to be the Riesz isomorphism associated
with the (equivalent) norm on $V$ given by
\beeq{norB}
  \nord{v}R=\dua{Rv,v}{V^*}V
   =\io |\nabla v|^2+\eta\iga v^2.
\end{equation}
At this point, we are able to state our first result,
related to existence and uniqueness of solutions (with
$\varepsilon=1$),
which slightly extends \cite[Thm.~1]{Ro}:
\bete\label{teoexist}
 Let\/ \eqref{hplambda}, \eqref{hpW}--\eqref{hpW2},
 \eqref{hpj}--\eqref{hpj2}, \eqref{hpteta0}, \eqref{hpchi0},
 \eqref{f1} hold. Moreover, assume either (Dirichlet conditions)
 \beeq{diri}\tag{Dirichlet}
   B:=-\Delta:V_0\to V_0^*, \quad\
   \calV:=V_0, \quad\
   \|\cdot\|_{\calV}:=\|\cdot\|_{V_0}, \quad\
   g:=f,
 \end{equation}
 or\/ \eqref{tetaga} and (Robin conditions)
 \beeq{rob}\tag{Robin}
   B:=R, \ \quad
   \calV:=V, \ \quad
   \|\cdot\|_{\calV}:=\|\cdot\|_{R}, \ \quad
   \dua{g,v}{\calV^*}\calV:=(f,v)+\eta\iga j'(\teta_\Gamma)\+v,
 \end{equation}
 the last relation holding for all $v\in\calV=V$, a.e. in $(0, \infty)$.
 Then, there exist a constant $c>0$, depending only on the data
 $\lambda,\, W,\,j,\,\teta_0,\,\chi_0$, and
 a\/ {\rm unique} pair $(\teta,\chi)$ such that if $u=j'(\teta)$
 \beal{regoteta}
   & \nor{\teta_t}{L^2(0,\infty;\calV^*)}
    +\nor{\teta}{L^\infty(0,\infty;H)}\le c,\\
  \label{regou}
   & \nor{u}{L^2(0,\infty;\calV)}\le c, \\
  \label{regoteta-2}
   & \nor{\teta-\teta_\infty}{L^2(0,\infty;\calV)}\le c,\\
  \label{regochi}
   & \nor{\chi_t}{L^2(0,\infty;H)}
    +\nor{\chi}{L^\infty(0,\infty;V)}\le c,
 \end{align}
 $\chi \in C^0([0,t];V) \cap L^2(0,t;H^2(\Omega))$
 for all $t>0$, and $\teta,\, \chi, \, u$ satisfy
 \beal{calore}
   & \teta_t+\lambda(\chi)_t+Bu=g \quext{in }\,\calV^*,\\
  \label{phase}
   & \chi_t+A\chi+W'(\chi)=\lambda'(\chi)u \quext{in }\,V^*,
 \end{align}
 a.e.~in $(0,\infty)$, as well as
 \beeq{iniz}
   \teta|_{t=0}=\teta_0,\qquad
   \chi|_{t=0}=\chi_0.
 \end{equation}
\ente
\noindent%
%
%
\beos\label{suregog}
 We point out that, in the \eqref{rob} case, \eqref{f1} and
 \eqref{tetaga} entail
 \beeq{gdaf}
  g\in L^2(0,\infty;\calV^*).
 \end{equation}
 This is of course true also in the \eqref{diri}
 case, in which we even have by assumption the better relation $g\in L^2(0,\infty;H)$.
\eddos
\noindent%
We shall not give the full proof of Theorem~\ref{teoexist}
since it is just a variant of the proof
given in \cite{Ro}. More precisely, the same argument of
\cite[Sec.~3]{Ro} can be used to obtain
existence of a solution to a
suitable regularization of the problem.
Concerning the a priori estimates which are required
to remove the approximation, some points
are technically different from \cite[Sec.~3]{Ro}
especially due to our choice(s) of boundary conditions.
Nevertheless, the required estimates will be
easily obtainable from the uniform bounds we shall
prove in the long-time analysis.
Finally, the compactness argument necessary to
pass to the limit and the proof of uniqueness
will be briefly sketched in the next~section.

Let us come now to our statement on
further regularity properties of solutions. To this aim,
we prepare an auxiliary result.
\bele\label{lemconv}
 Let X be a Banach space and let
 $\gamma\in L^2(0,\infty;X)$ satisfy, for
 some $p\in[1,\infty]$,
 \beeq{g}
   \gamma\in W^{1,p}(0,t;X) \quad \perogni t \in (0,\infty), \qquad
    \sup_{t\in[0,\infty)} \nor{\gamma_t}{L^p(t,t+1;X)}<\infty.
 \end{equation}
Then, it turns out that $\gamma\in L^\infty(0,\infty;X)$.
Moreover, if $p>1$, the strong convergence $\gamma(t)\to 0$ holds
in $X$, as $t\nearrow \infty$.
\enle
\noindent%
{\bf Proof.}~~%
First of all, the $L^\infty$ bound
can be shown using the Fundamental Theorem
of Calculus, by a simple contradiction argument.
Concerning the convergence to $0$, let $p>1$ and
$\{t_n\}$ be
an arbitrary diverging sequence of times. Setting
$$
  \gamma_n(t):=\gamma(t_n+t), \qquad t\in[0,1],
$$
it is clear that $\gamma_n\to 0$ strongly in $L^2(0,1;X)$. Thus,
there is a subsequence, not relabelled, such that
$\gamma_n(t)\to0$ strongly in $X$ for a.e.~$t\in(0,1)$.
In particular, we can find a sequence $\{\delta_k\}\subset(0,1)$,
with $\delta_k\searrow0$ and such that, in $X$,
$$
  \lim_{n\to\infty} \gamma_n(\delta_k)=0
   \qquad\perogni k\in\NN.
$$
We then deduce from \eqref{g}
$$
  \|\gamma(t_n)\|_X
   \le \|\gamma(t_n)-\gamma(t_n+\delta_k)\|_X
        +\|\gamma(t_n+\delta_k)\|_X
   \le c\delta_k^{1/p^*}+\|\gamma_n(\delta_k)\|_X,
$$
where $p^*$ is the conjugate exponent to $p$,
and the latter quantity can be made arbitrarily small for
$n$ large enough. Due to the arbitrariness of $\{t_n\}$, this
shows that $\gamma(t)\to 0$ in $X$ as $t\nearrow\infty$.%
\dimbox
\beos\label{ossconv}
 We point out that the second part of the
 statement above is false (even for
 $X=\RR$) if one takes $p=1$.
 Indeed, set, for $n\ge 2$, $v_n:\RR\to\RR$ given by
 $$
   v_n(t):=n^2\calX_{(-1/n^2,0)}-n^2\calX_{(0,1/n^2)},
 $$
 $\calX$ denoting the characteristic function,
 and define
 $$
  \gamma(t):=\sum_{n=2}^\infty v_n(t-n),\quad t\in \RR,
   \qquad g(t):=\int_0^t \gamma(s)\,\dis
 $$
 Then, it is clear that $g$ stays in
 $L^2(0,\infty)\cap L^\infty(0,\infty)$ and
 satisfies \eqref{g} (with $X=\RR$).
 However, $g(t)$ does not tend to $0$ for $t\nearrow\infty$.
\eddos
\bete\label{teorego}
 Let\/ \eqref{hplambda}, \eqref{hpW}--\eqref{hpW2},
 \eqref{hpj}--\eqref{hpj2}, \eqref{hpteta0},
 \eqref{hpchi0}, \eqref{f1} and either\/ \eqref{diri}, or\/
 \eqref{tetaga} and\/ \eqref{rob}, hold.
 Let also \eqref{g} hold for $\gamma=g$, $X=\calV^*$
 and some $p\in[1,\infty]$. Then, for all
 $s>0$ there exists a constant $c>0$,
 depending on $\lambda,\, W,\, j,\, \teta_0,\, \chi_0$ and $s$,
 such that
 \beal{regoteta2}
   & \sup_{t\ge s}\nor{\teta_t}{L^2(t,t+1;H)}
    +\nor{\teta}{L^\infty(s,\infty;\calV)}
    \le c,\\
  \label{regou2}
   & \nor{u}{L^\infty(s,\infty;\calV)}\le c,\ \hbox{ with }\, u=j'(\teta),\\[1mm]
  \label{regochi2}
   & \nor{\chi_t}{L^{\infty}(s,\infty;H)}
    +\sup_{t\ge s}\nor{\chi_t}{L^2(t,t+1;V)}
    +\nor{\chi}{L^\infty(s,\infty;H^2(\Omega))}\le c,\\
  \label{regoW2}
   & \nor{W'(\chi)}{L^\infty(s,\infty;H)}\le c.
 \end{align}
 If, additionally,
 \beeq{newregog}
   g_t\in L^q(0,\infty;\calV^*)
 \end{equation}
 for some $q\in[1,2]$, then we also have
  \beeq{regoadd-1}
   \nor{\teta_t}{L^2(s,\infty;H)}
    +\nor{\chi_t}{L^2(s,\infty;V)}
    \le c.
 \end{equation}
\ente
\noindent%
In particular, the above theorem provides
a priori estimates
which have a uniform character for large times.
Of course, this uniformity cannot be proved for the
$L^2$ in time norm of $\chi$ in \eqref{regochi}.

Let us now move to the study of long-time behavior, starting
from existence of a nonempty $\omega$-limit set.
\bete\label{teoomega1}
 Let\/ \eqref{hplambda}, \eqref{hpW}--\eqref{hpW2},
 \eqref{hpj}--\eqref{hpj2}, \eqref{hpteta0},
 \eqref{hpchi0}, \eqref{f1} and either\/ \eqref{diri}, or\/
 \eqref{tetaga} and\/ \eqref{rob}, hold.
 Let also \eqref{g} hold for $\gamma=g$, $X=\calV^*$
 and some $p\in[1,\infty]$.
 Moreover, let us assume that, either
 \beeq{gl2}
   g_t\in L^2(0,\infty;\calV^*)
 \end{equation}
 or
 \beeq{j12}
   \esiste c>0,\, \alpha\le 3:~~
   |j''(r)|\le c\big(1+ |j'(r)|^\alpha \big)
   \quad\perogni r\in J.
 \end{equation}
 Then, we have that, as $t\nearrow\infty$,
 \beeq{utetainfi}
   u(t)\to 0,~~\teta(t)\to\teta_\infty
    \quext{weakly in }\,\calV~~
    \text{and strongly in }\,H.
 \end{equation}
 Moreover, any diverging sequence $\{t_n\}\subset(0,\infty)$
 admits a subsequence, not relabelled, such that
 \beeq{chiinfi}
   \chi(t_n)\to\chi_\infty \quext{weakly in }\,H^2(\Omega)~~
    \text{and strongly in }\,V\cap C(\overline\Omega),
 \end{equation}
 where $\chi_\infty$ is a solution of the stationary
 problem
 \beeq{prostaz}
   A\chi_\infty+W'(\chi_\infty)=0  \quext{in }\,V^*.
 \end{equation}
\ente
\noindent%
\beos\label{reinf}
 Assumption \eqref{g} on $g$ of course
 reinforces the convergence properties of
 $g$ to $0$ (cf.~\eqref{f1} and,
 in the \eqref{rob} case, \eqref{tetaga}).
\eddos
\noindent%
\beos\label{qualij}
 The property \eqref{j12} is actually not very strong. For instance,
 in the situation when $J=\RR$, then \eqref{j12} is fulfilled provided
 $j$ has a polynomial, or even exponential, growth at infinity.
 Moreover, \eqref{j12} also holds for the laws mentioned
 in the Introduction (corresponding to combinations of the
 Caginalp and Penrose-Fife models).  A case in which
 \eqref{j12} does not hold (but \eqref{hpj}--\eqref{hpj2} do hold)
 is given by $j'(r)=r+\log (1 + r/\tau_c)$, due to the singular
 behavior of $j'$ in proximity of $- \tau_c $.
\eddos
\noindent
Our last result characterizes the $\omega$-limit as a singleton
in case the potential $W$ is {\sl analytic}.
To introduce it we need some preliminaries. First of all,
let us notice that \eqref{hpW}--\eqref{hpW2} entail by
simple maximum principle arguments (see \cite[Lemma~3.1]{AF})
that there exists a closed interval $I_1\subset I_0$ such
that any solution $\chi_\infty$ to \eqref{prostaz}
satisfies
\beeq{maxprin}
   \chi_\infty\in W^{2,q}(\Omega)~~\perogni q\in [1,\infty),\qquad
  \chi_\infty(x)\in I_1\quext{for all }\,x\in\Omega.
\end{equation}
Moreover, if we set, for $v\in V$,
\beeq{enloj}
  E(v):=\io\frac12|\nabla v(x)|^2+W(v(x))\,\dix
\end{equation}
(which might be $+\infty$ if $W(v)$ is not summable),
there holds the following form of the {\sl Simon-\L ojasiewicz}\/
inequality, which is a reformulation of
\cite[Prop.~4.4]{AF}
(see also \cite[Prop.~4.2]{AFIR}):
\bete\label{teoloj}
 Let \eqref{hpW}--\eqref{hpW2} hold and let
 \beeq{anal}
   W\,\text{ be real analytic on }\,I_0.
 \end{equation}
 Let $\chi_\infty$ be a solution
 to \eqref{prostaz}. Then,
 there exist constants $c_\ell,\epsilon>0$, $\zeta\in(0,1/2)$,
 such that
 \beeq{disloj}
   |E(v)-E(\chi_\infty)|^{1-\zeta} \le c_\ell
    \|Av+W'(v)\|_{V^*}
 \end{equation}
 for all $v\in V$ such that
 \beeq{disloj2}
    \|v-\chi_\infty\|_{V\cap C(\overline\Omega)}\le\epsilon.
 \end{equation}
\ente
\beos\label{ossbarriere}
 The above statement is given in a sligthly different
 fashion with respect to \cite[Prop.~4.4]{AF} since this
 version seems to be more suitable for our specific problem.
 However, we point out that our hypotheses
 entail those of \cite[Prop.~4.4]{AF}. Actually,
 by \eqref{maxprin}, it
 is clear that, as we possibly take a smaller $\epsilon$
 in condition \eqref{disloj2}, then any $v$ fulfilling
 \eqref{disloj2} also satisfies
 \cite[(4.5)]{AF}. Indeed, $\epsilon$ can be taken so
 small that $v$ ranges into $I_0$, where Lipschitz
 continuity of $W$ holds (recall that
 $\overline{I_0}\subset I$).
\eddos
\beos\label{ossbarriere2}
 Notice also that in this way we actually get rid of the
 possibly singular character of $W$ at the boundary of
 $I$. Indeed, it is not excluded that there is a
 transient dynamics where $W'(\chi(t))$ may be
 unbounded. However, thanks to \eqref{maxprin} and
 the precompactness of the trajectory in $C(\overline\Omega)$,
 for sufficiently large times $\chi(t)$ takes values
 into $I_0$, where $W$ is bounded and analytic.
 As in \cite{GPS1}, the key condition
 ensuring this property  is \eqref{hpW2},
 which essentially states that the leftmost
 and rightmost minima of $W$ are {\sl interior}\/ to
 its domain. The gradient flow structure of the
 system entails that the solution eventually moves
 away from these minima.
\eddos
\noindent%
Here is, finally, our convergence result:
\bete\label{teoomega2}
 Let the hypotheses of\/ {\rm Theorem~\ref{teoomega1} hold.}
 Furthermore, assume\/ \eqref{anal} and
 \begin{equation}\label{effe2}
   \sup_{t\ge 0} t^{1+\delta}\,\int_t^{\infty} \Vert g(s)\Vert_{\calV^*}^2 \,ds < \infty,
 \end{equation}
 for some $\delta>0$. Then, as $(\teta_0,\chi_0)$
 are initial data satisfying
 \eqref{hpteta0}--\eqref{hpchi0}, the $\omega$-limit
 of the corrisponding trajectory $(\teta,\chi)$ of system\/
 \eqref{calore}--\eqref{phase} consists of a\/ {\rm unique} pair
 $(\teta_\infty,\chi_\infty)$, where $\teta_\infty$ is given by\/
 \eqref{hpj2} and $\chi_\infty\in V$ is a solution to\/
 \eqref{prostaz}. Moreover, as $t\nearrow\infty$,
 \eqref{utetainfi} holds together with
 \beeq{chiinfi2}
   \chi(t)\to\chi_\infty \quext{weakly in }\,H^2(\Omega)~~
    \text{and strongly in }\,V\cap C(\overline\Omega).
 \end{equation}
 More precisely, if
 \begin{equation}\label{condelta1}
   \delta > \frac{2\zeta}{1-2\zeta},
 \end{equation}
 where $\zeta$ is as in \eqref{disloj},
 then one can find $t^*>0$ and a positive constant
 $c_*$ such that
 \begin{equation}\label{rate1}
   |\chi(t) - \chi_\infty| \le c_* t^{-\frac{\zeta}{1-2\zeta}},
    \qquad \forall\,t\ge t^*.
 \end{equation}
 Otherwise, one can find $\zeta_0\in(0,\zeta)$ so that
 \begin{equation}\label{condelta2}
   \delta > \frac{2\zeta_0}{1-2\zeta_0},
 \end{equation}
 a time $t^{**}>0$ and a positive constant $c_{**}$ such that
 \begin{equation}\label{rate2}
   |\chi(t) - \chi_\infty| \le c_{**} t^{-\frac{\zeta_0}{1-2\zeta_0}},
    \qquad \forall\,t\ge t^{**}.
 \end{equation}
\ente
\noindent
\beos\label{convinV}
 Let us notice that \eqref{rate1}, \eqref{rate2} give the convergence
 rate of $\chi(t)$ to $\chi_\infty$ with respect to the norm
 of $H$. Of course, an estimate of the rate of convergence
 in the norm of $V$ could be obtained from the uniform bound
 corresponding to the last of \eqref{regochi2} and
 interpolation.
\eddos


\section{A priori estimates and well posedness}
\label{secunif}

Let us first sketch an approximated version of system
\eqref{calore}--\eqref{phase} along the lines of \cite{Ro},
\cite[Sec.~3]{FSc}, to which we refer for more details.
Namely, let us assume that $j$ and $W$ have been
replaced in \eqref{calore}--\eqref{phase}
by regularized functions $j_n$ and $W_n$
defined on the whole real line
and such that
\beeq{Mosco}
  j_n,\, W_n + \frac{\kappa}2 \,\Id^2 \to j,\, W + \frac{\kappa}2 \, \Id^2
   \quext{in the sense of Mosco}
\end{equation}
(see, e.g., \cite{At} for the definition of Mosco convergence
and for the related notion of $G$-convergence of graphs).
Moreover, we can assume \cite[Sec.~3]{FSc} that, for all
$n\in\NN$,
\bega{Wn}
  W_n''(r)\ge -\kappa~~\perogni r\in \RR, \qquad
   \frac{W_n'(r)}r \ge \mu/2~~\perogni r\in \RR\setminus I_0,\\
 \label{jn}
   j_n''(r)\ge \sigma/2~~\perogni r\in\RR,
\end{gather}
where $\mu$, $I_0$, and $\sigma$ are as in \eqref{hpW},
\eqref{hpW2}, \eqref{hpj} (cf. \cite[Sec.~3]{FSc} for an
example of a possible regularizing sequence $\{W_n\}$).
Then, we consider a family $\{(\tetan,\chin)\}_{n\in\NN}$
of (possibly local in time) solutions to
the regularized problem specified by the subscript $n$.
Existence of these solutions can be shown
proceeding as in \cite[Subsec.~3.1]{Ro}
and it turns out that $(\tetan,\chin)$, as well
as $u_n:=j_n'(\tetan)$ and $W_n'(\chin)$, are regular
enough to give a rigorous meaning to the
forthcoming computations. Actually, we shall now deduce
some a priori estimates with the aim of
taking the limit as $n\nearrow\infty$.
In this procedure, \eqref{calore} and  \eqref{phase}
will be implicitly considered in their $n$-approximated
form. Moreover, $c$ will denote a positive constant,
whose value is allowed to vary even inside
one single row, but $c$ may depend only on
$\lambda,j,W,\teta_0,\chi_0$ (and
neither on $n$ nor on $t$).
When we need to fix the
value of some specific $c$, we shall use the
notation $c_i$, $i\ge 0$.
The symbols $C$, $C_i$ will denote constants
which, instead, can explicitly depend on $t\in(0,\infty)$,
but do not explode as $t\searrow 0$.
For simplicity, 
we shall proceed as if the solutions were defined
for all times $t\ge 0$.
Indeed, although this might be not true at the
approximating level, it will certainly hold at the limit in
view of the uniform in time
character of the estimates and of
standard extension arguments.

\vspace{2mm}

\noindent%
{\bf Energy estimate.}~~%
Test \eqref{calore} by
$u_n$ in the duality between $\calV^*$
and $\calV$ and sum the result to
\eqref{phase} tested by $\chint$ in the
duality of $V^*$ and $V$.
The smoothness properties assumed on the approximating
solutions, standard integration by parts formulas,
and the cancellation of a couple of opposite
terms then give
\beeq{conto11}
  \ddt\calE_n(\tetan,\chin)
   +|\chint|^2
   +\frac12\nord{u_n}{\calV}
   \le\frac12\|g\|_{\calV^*}^2,
\end{equation}
where we also used the Young 
inequality 
to split the
duality product of $g$ and $u_n$ resulting from the
\rhs\ of \eqref{calore}. Also, accordingly
with \eqref{energyin}, we have set
\beeq{energyn}
  \calE_n(\teta,\chi):=\io\Big(
   \frac12|\nabla\chi|^2
   +W_n(\chi)
   +j_n(\teta)\Big).
\end{equation}
Owing now to \eqref{hpteta0}, \eqref{hpchi0}, and
\eqref{g}, we can integrate
\eqref{conto11} over $(0,t)$ for arbitrary
$t>0$ and deduce
\beal{st11}
  & \nor{u_n}{L^2(0,t;\calV)}
  +\nor{j_n(\tetan)}{L^\infty(0,t;L^1(\Omega))}\le c,\\
 \label{st12}
  & \nor{\chint}{L^2(0,t;H)}
  +\nor{\chin}{L^\infty(0,t;V)}
  +\nor{W_n(\chin)}{L^\infty(0,t;L^1(\Omega))}\le c.
\end{align}
Let us note that, to obtain the second of \eqref{st12}, we
used that, by \eqref{Wn},
\beeq{encoerc}
  \esiste c,\, c_0>0:\quad\io\Big(
   \frac12|\nabla v|^2
   +W_n(v)\Big)\ge c\|v\|_V^2-c_0
   \quad\perogni v\in V.
\end{equation}

\vspace{2mm}

\noindent%
{\bf Second estimate.}~~%
Let us now test \eqref{calore} by $\tetan-\teta_\infty$.
Owing to the monotonicity of $j_n'$ and, more precisely,
to \eqref{jn}, we can see that, both in
the \eqref{diri} and in the \eqref{rob}~case,
\beeq{conscoercjn}
  \dua{B u_n,\tetan-\teta_\infty}{\calV^*}\calV
   \ge 2c_1\|\tetan-\teta_\infty\|_{\calV}^2,
\end{equation}
for some $c_1>0$ depending in particular on $\sigma$. Thus,
splitting the term depending on $g$ by the Young inequality,
we get
\beeq{conto31}
  \frac12\ddt|\tetan-\teta_\infty|^2
   +c_1\|\tetan-\teta_\infty\|_{\calV}^2
   \le c\|g\|_{\calV^*}^2
   -\dua{\lambda'(\chin)\chint,\tetan-\teta_\infty}{\calV^*}\calV,
\end{equation}
and the latter term is readily estimated as follows
\bealo
  {}-\dua{\lambda'(\chin)\chint,\tetan-\teta_\infty}{\calV^*}\calV
   & = {}-(\lambda'(\chin)\chint,\tetan-\teta_\infty)\\
 \no
  & \le {}c\big(1+\nor{\chin}{L^4(\Omega)}\big)
    |\chint|\+\|\tetan-\teta_\infty\|_{L^4(\Omega)}\\
 \label{conto32}
  & \le {}c\big(1+\nord{\chin}{V}\big)|\chint|^2
   +\frac{c_1}2\|\tetan-\teta_\infty\|_{\calV}^2,
\end{align}
where we used \eqref{hplambda} and the continuous embeddings
$V,\, \calV\subset L^4(\Omega)$.
Hence, recalling \eqref{st12}
and using \eqref{gdaf}, we infer
\beeq{st31}
  \nor{\tetan}{L^\infty(0,t;H)}
   +\nor{\tetan-\teta_\infty}{L^2(0,t;\calV)}
   \le c.
\end{equation}
Moreover, arguing as
in \eqref{conto32}, we obtain that the term
$\lambda(\chin)_t=\lambda'(\chin)\chi_{n,t}$
in \eqref{calore} is uniformly
controlled in $L^2(0,T; \calV^*)$, so that
\eqref{st11}, \eqref{gdaf}
and the equality $\teta_{n,t}= - \lambda'(\chin)\chi_{n,t} -Bu_n +g$
yield
\beeq{st42}
  \nor{\teta_{n,t}}{L^2(0,T;\calV^*)}\le c .
\end{equation}

\vspace{2mm}

\noindent%
{\bf Third estimate.}~~%
Thanks to \eqref{hplambda}, the continuous embeddings $V,\, \calV\subset L^4(\Omega)$,
and \eqref{st11}--\eqref{st12}, we infer that
\beeq{st21}
  \nor{\lambda'(\chin)u_n}{L^2(0,t;H)}
   \le c\big(1+\nor{\chin}{L^\infty(0,t;L^4(\Omega))}\big)
   \nor{u_n}{L^2(0,t;L^4(\Omega))}\le c.
\end{equation}
Then, we can test \eqref{phase} by $A\chin$
and integrate over $\Omega\times(0,t)$. Using
the first inequality of \eqref{Wn}, we have that
\beeq{nonpar}
  \itt\big(A\chin,W_n'(\chin)\big)
   \ge -\kappa\|\chin\|_{L^2(0,t;V)}^2,
\end{equation}
which can be controlled thanks to \eqref{st12},
but only on bounded time intervals. Therefore, using
also \eqref{hpchi0} it is not difficult to obtain
\beeq{st22}
  \nor{\chin}{L^2(0,t;H^2(\Omega))}
   +\nor{W_n'(\chin)}{L^2(0,t;H)}
   \le C,
\end{equation}
where the second bound comes from a further comparison
of terms in \eqref{phase}.

\vspace{2mm}

\noindent%
{\bf Limit as $n\to +\infty$ and existence.}~~%
Standard compactness tools enable us to pass to the limit in the
$n$-approximated versions of \eqref{calore}--\eqref{phase}. Indeed,
thanks to estimates \eqref{st11}--\eqref{st12}, \eqref{st31}--\eqref{st42},
and \eqref{st22}, there exist
four limit functions $\teta, u, \chi, v,$ defined from $(0,\infty)$
to $H$ (at least), and a suitable subsequence of $n$ (not relabeled)
such that the corresponding subsequences $\{\tetan\}$, $\{u_n\}$, $\{\chin\}$,
$\{W'_n(\chin)\}$ fulfill
\beal{conv1}
  & \teta_n \to \teta \ \hbox{ weakly star in }  \ H^1(0,t;\calV^*)
\cap L^\infty(0;t; H) \cap L^2 (0,t;\calV) ,\\
\label{conv2}
  & u_n \to u \ \hbox{ weakly in }  L^2 (0,t;\calV),\\
\label{conv3}
  & \chin \to \chi \ \hbox{ weakly star in }  \ H^1(0,t;H)
\cap L^\infty(0;t; V) \cap L^2 (0,t;H^2(\Omega)) ,\\
\label{conv4}
  & W'_n(\chin) \to v \ \hbox{ weakly in }  L^2 (0,t;H)
\end{align}
as $n\nearrow \infty$, for all $t>0$. We note at once that the bounds
in \eqref{regoteta}--\eqref{regochi} are certainly satisfied by the limit
functions: in fact, it suffices to take the $\liminf$ in estimates
\eqref{st11}--\eqref{st12},
\eqref{st31}--\eqref{st42} which are uniform with respect to $t$.
Next, \eqref{conv1}, \eqref{conv3},
the Ascoli theorem and the Aubin compactness lemma (see, e.g.,
\cite[Cor.~4, Sec.~8]{Si}) enable us to deduce that
\beal{convf1}
  & \teta_n \to \teta \ \hbox{ strongly in } \
   C^0([0;t];\calV^* ) \cap L^2 (0,t;H) ,\\
\label{convf2}
  & \chin \to \chi \ \hbox{ strongly in }  \
   C^0([0;t]; H) \cap L^2 (0,t;\calV ) ,
\end{align}
whence $ \lambda'(\chin) \to \lambda'(\chi)$ strongly in $C^0([0;t]; H)$
due to the Lipschitz continuity of $\lambda'$. Then,
\eqref{conv2}--\eqref{conv3} and the continuous embedding
$V\subset L^6(\Omega)$ imply
\beal{con-la1}
  & \lambda'(\chin)u_n \to \lambda'(\chi) u \ \hbox{ weakly in }  \
    L^2 (0,t;L^{3/2} (\Omega) ) ,\\
\label{con-la2}
  & \lambda(\chin)_t=\lambda'(\chin)\chi_{n,t} \to
  \lambda(\chi)_t=\lambda'(\chi)\chi_{t} \ \hbox{ weakly in } \
  L^2 (0,t;L^1(\Omega)) .
\end{align}
At this point, we can pass to the limit in the
$n$-approximated versions of \eqref{calore}--\eqref{phase}
to obtain \eqref{calore} (which a fortiori holds in $\calV^*$)
and
\beeq{phase-v}
   \chi_t+A\chi+ v=\lambda'(\chi)u \quext{in }\,V^*,
\end{equation}
a.e. in $(0,\infty)$. Initial conditions \eqref{iniz} follow
easily from \eqref{convf1}--\eqref{convf2}. Then, in order to
conclude the existence proof, it remains to identify
functions $u$ and $v$, that is, to check that
\beeq{ident}
   u = j'(\teta), \quad v= W'(\chi)  \quad \text{ a.e. in }\, \Omega\times (0,\infty).
\end{equation}
However, due to the Mosco convergences in \eqref{Mosco},
it turns out that (cf., e.g., \cite[Thm.~3.66]{At})
the subdifferential operators \ $j'_n, \, W'_n + \kappa \Id$  \
$G$-converge to \ $j', \, W' + \kappa \Id$, \ as well as their extensions
to $L^2(\Omega\times (0,t))$, for all $t>0$. Then, we can apply the basic
properties of $G$-convergence (see, e.g. \cite[Prop.~1.1, p.~42]{Ba}) stating
that if $ u_n= j'_n (\teta_n) \to u$ and $\teta_n \to \teta$ weakly in
$L^2(\Omega\times (0,t))$ and $\limsup \int\!\!\int_{\Omega\times (0,t)} u_n\teta_n \leq
\int\!\!\int_{\Omega\times (0,t)} u\teta$, then $u = j'(\teta)$. But, in our case this
follows easily from \eqref{conv2} and \eqref{convf1}. The other identification in
\eqref{ident} is a bit longer, since we first check that $v + \kappa \chi \in
(W' + \kappa \Id)(\chi) $ with the help of \eqref{conv4} and
\eqref{convf2} (recall that $W' + \kappa \Id$ is monotone
by \eqref{hpW}), then extract the information $ v= W'(\chi)$.
Hence, we conclude for the validity of~\eqref{ident}.

\vspace{2mm}

\noindent%
{\bf Uniqueness.}~~%
Assume, by contradiction, that there exist two solutions $(\teta_i,\chi_i)$, $i=1,2$,
to problem \eqref{calore}--\eqref{iniz}, let $u_i= j'(\teta_i)$, $i=1,2$, and set
temporarily $\teta:=\teta_1 - \teta_2$,  $\chi:=\chi_1- \chi_2$, $u:=u_1 - u_2$.
Now, we can take the difference of \eqref{calore}
written for the solutions corresponding to $i=1,2$,
integrate it in time from $0$ to $t>0$,
and test by $u(t)$. At the same time, we test the difference of
\eqref{phase} by $\chi(t) $ and sum the result to the previous relation.
With the help of \eqref{hpj} and \eqref{hpW} it is straightforward to deduce
\bealo
  & \sigma |\teta(t)|^2 + \unmezzo \ddt
   \left(\big\| {\textstyle \int_0^t u(s) ds }\big\|^2_{\calV}
    +|\chi(t)|^2 \right) +|\nabla\chi (t)|^2\\
  & 
   \le {} \kappa|\chi(t) |^2
   + \io \big((\lambda'(\chi_1)u_1 - \lambda'(\chi_2) u_2 )\chi - (\lambda(\chi_1)-\lambda(\chi_2))u \big)(t).
 \label{dis-uni}
\end{align}
To estimate the last term we use a remark from Kenmochi (see \cite{Ken} and, e.g.,
\cite[Lemma~3.2]{KK}): in fact, by the Taylor expansion and \eqref{hplambda} we have
\begin{align*}
&
(\lambda'(\chi_1)u_1 - \lambda'(\chi_2) u_2 )\chi - (\lambda(\chi_1)-\lambda(\chi_2))u \\
&{}= u_1 (\lambda(\chi_2) - \lambda(\chi_1) - \lambda'(\chi_1) (\chi_2 - \chi_1) )
+   u_2 (\lambda(\chi_1) - \lambda(\chi_2) - \lambda'(\chi_2) (\chi_1 - \chi_2) ) \\
&{}\leq \frac{\Lambda}2 (|u_1| + |u_2| )|\chi|^2.
\end{align*}
Therefore, the \lhs\ of \eqref{dis-uni} can be handled using the H\"older inequality
and the continuous embeddings $V,\, \calV\subset L^4(\Omega)$ in order to obtain
\bealo
  & \sigma \int_0^t |\teta(s)|^2 ds  + \unmezzo
   \left(\big\| {\textstyle \int_0^t u(s) ds }\big\|^2_{\calV}
    +|\chi(t)|^2 \right) + \int_0^t |\nabla\chi (s)|^2 ds\\
  & 
   \le {} c  \int_0^t \big( 1+ \|u_1(s) \|_{\calV}^2
             + \|u_2 (s) \|_{\calV}^2 \big) |\chi(s) |^2 ds
   +\unmezzo \int_0^t \|\chi (s) \|_V^2 ds .
 \label{disII}
\end{align}
Finally, in view of the regularity \eqref{regou} for $u_1, \, u_2 $,
the uniqueness property follows easily from the Gronwall lemma.
This completes the proof of  Theorem~\ref{teoexist}.


\vspace{2mm}

\noindent%
{\bf Proof of Theorem~\ref{teorego}.}~~%
Let us test \eqref{calore} by
$u_t$ in the duality between $\calV^*$
and $\calV$. Then, differentiate in time
\eqref{phase}, multiply the result by $\chi_t$,
and integrate over $\Omega$. Summing together
the obtained relations, noting that
a couple of terms cancel out, and using
\eqref{hpW}, we infer
\bealo
  & \ddt\Big(
   \frac12\|u\|_{\calV}^2
   +\frac12|\chi_t|^2
   -\dua{g,u}{\calV^*}{\calV}
        \Big)
   +(u_t,\teta_t)
   +|\nabla\chi_t|^2\\
 \label{conto41}
  & \mbox{}~~~~~~~~~
   \le-\dua{g_t,u}{\calV^*}{\calV}
   +\kappa|\chi_t|^2
   +\io\lambda''(\chi)\chi_t^2u.
\end{align}
Let us note that the computation above is just
formal in the regularity setting of Theorem~\ref{teoexist}.
However, the procedure might be made rigorous by working on
the $n$-regularization sketched before
and then passing to the limit.
We omit the details just for brevity.
By \eqref{hpj}, we have
\beeq{conto42}
  (u_t,\teta_t)
    \ge\sigma|\teta_t|^2 .
\end{equation}
Moreover, using \eqref{hplambda}
and once more the continuous embeddings
$V,\, \calV\subset L^4(\Omega)$,
we obtain
\beeq{conto43}
  \io\lambda''(\chi)\chi_t^2u
   \le\frac12\|\chi_t\|_V^2
    +c\|u\|^2_\calV|\chi_t|^2.
\end{equation}
Finally, we note that
\beeq{conto44}
  -\dua{g_t,u}{\calV^*}{\calV}
   \le\big\|g_t\|_{\calV^*}\|u\|_{\calV}.
\end{equation}
Next, we set
\beeq{defF}
 \calF:=
   \frac12\|u\|_{\calV}^2
   +\frac12|\chi_t|^2
   -\dua{g,u}{\calV^*}{\calV}
  \end{equation}
 and observe that,
since $g\in L^\infty(0,\infty;\calV^*)$ thanks to
Lemma~\ref{lemconv}, there exist constants $c_2 , c_3>0$
depending only on $g$ and such that
\beeq{stimaF}
   \frac14\|u\|_{\calV}^2
    +\frac12|\chi_t|^2
   \le\calF+c_2\le
  \|u\|_{\calV}^2
  +\frac12|\chi_t|^2+c_3
\end{equation}
for all $t>0$.
Thus, setting $\calY:=\calF+c_2$
and using \eqref{conto42}--\eqref{conto44},
\eqref{conto41} becomes
\bealo
  \ddt \calY
  +\sigma|\teta_t|^2
  +\frac12|\nabla\chi_t|^2
  & \le c|\chi_t|^2\big(1+\|u\|_{\calV}^2\big)
   +\|u\|_{\calV}\|g_t\|_{\calV^*}\\
 \no
  & \le c \calY (1+\calY)
   +\big(1+\|u\|_{\calV}^2\big)\big(\|g_t\|_{\calV^*}^p+1\big)\\[1mm]
 \label{conto45}
  & \le c(1+\calY)\big(1+\calY+\|g_t\|_{\calV^*}^p\big).
\end{align}
Then, noting that, since  by the
estimates \eqref{regou}--\eqref{regochi},
\beeq{einT1}
  \sup_{t\ge 0}\int_t^{t+1}
   \calY(s)\,\dis<\infty,
\end{equation}
we can apply to $\calY$ the {\sl uniform}\/
Gronwall Lemma (see, e.g., \cite[Lemma III.1.1]{Te}),
which yields \eqref{regou2} and the first bound
in \eqref{regochi2}. Next, integrating \eqref{conto45}
in time over $(t,t+1)$ for $t$ greater than or equal to
a given $s>0$, we get the first of \eqref{regoteta2}
and the second of \eqref{regochi2}.
Observing that, by \eqref{hpj}--\eqref{hpj2} (and the
Poincar\'e inequality in the \eqref{diri} case), there exists
$c>0$ such that
\beeq{nortetau}
  \|u\|_\calV\ge c\|\teta-\teta_\infty\|_\calV,
\end{equation}
we get the second of \eqref{regoteta2}
from \eqref{regou2}.
Finally, by \eqref{hplambda} we deduce
\beeq{st2membro}
  \|\lambda'(\chi)u\|_{L^\infty(s,T;H)}
   \le c\big(1+\|\chi\|_{L^\infty(s,T;V)}\big)
   \|u\|_{L^\infty(s,T;\calV)}\le c.
\end{equation}
Thus, using the first of \eqref{regochi2} and \eqref{hpW}
and viewing \eqref{phase} as a time dependent family of elliptic
equations with monotone (up to a linear perturbation) nonlinearities,
the last of \eqref{regochi2} follows
from the standard elliptic regularity theory.
This, clearly, also gives \eqref{regoW2}. To conclude, we have
to prove \eqref{regoadd-1} in the case when \eqref{newregog}
holds. Coming back to the first row
of \eqref{conto45} it then suffices to note that,
for a.e.~$t\ge s$,
the \rhs\ is
\beeq{new-1}
  \le c|\chi_t|^2\big(1+\|u\|_{\calV}^2\big)
   +c\|u\|_{\calV}^{q^*}+c\|g_t\|_{\calV^*}^q
  \le c|\chi_t|^2+c\|g_t\|_{\calV^*}^q+c\|u\|_{\calV}^{q^*},
\end{equation}
where it is intended that, if $q\in(1,2]$,
then $q^*\in[2,\infty)$ denotes the conjugate exponent
to $q$; otherwise, i.e., if $q=1$, the latter term
on the \rhs\ has to be omitted. Relation
\eqref{regoadd-1} follows now by integrating \eqref{conto45}
over $(s,\infty)$ and using \eqref{new-1}.
Indeed, the terms on the \rhs\
of \eqref{new-1} are controlled, respectively, by
the first of \eqref{regochi}, by
\eqref{regou}, \eqref{regou2} and interpolation,
and by \eqref{newregog}. The
proof of the Theorem is now complete.\dimbox


\section{Study of the ${\boldsymbol{\omega}}$-limit}
\label{secomegalim}

{\bf Proof of Theorem~\ref{teoomega1}.}~~
Let us first notice that the convergence in \eqref{chiinfi}
to some $\chi_\infty \in H^2(\Omega)$ is an immediate consequence
of \eqref{regochi2} since $\chi$ is weakly continuous from $[s,\infty)$ to
$H^2(\Omega)$, whence the estimate $\|\chi(t)\|_{H^2(\Omega)} \le c$ holds
true for {\sl all}\/ $t\in[s,\infty)$.

Analogously,  the convergence $\teta(t)\to\teta_\infty$ in \eqref{utetainfi}
follows from \eqref{regoteta2} once one sees that  $\teta(t)\to\teta_\infty$
strongly in $\calV^*$ by \eqref{regoteta}, \eqref{regoteta-2},
and Lemma~\ref{lemconv}: indeed, due to the identification of the limit $\teta_\infty$,
it turns out that the entire family $\teta (t)$ converges.

Next, we show the convergence of $u(t)$ in
\eqref{utetainfi}. With this purpose, let us observe that, if \eqref{j12} holds,
then
\bealo
  \|u_t(t)\|_{L^1(\Omega)}
   & \le |j''(\teta(t))||\teta_t(t)|
   \le c\big(1+\|u(t)\|_{L^{2\alpha}(\Omega)}^\alpha\big)|\teta_t(t)|\\
 \no
  & \le c\big(1+\|u\|_{L^\infty(1,\infty;\calV)}^\alpha\big)
    |\teta_t(t)|
    \le c|\teta_t(t)|.
\end{align}
for a.e.~$t\in (1,\infty)$, say. Thus, by \eqref{regou} and \eqref{regoteta2}
it is clear that $\gamma=u$ fulfills the assumptions of
Lemma~\ref{lemconv} with $p=2$ and $X=L^1(\Omega)$ and
consequently $u(t)\to 0$ in $L^1(\Omega)$ as $t\nearrow\infty$.
Then, to show that \eqref{utetainfi} holds
(i.e., $u(t)$ weakly converges in $\calV$),
it is now enough to  point out the bound in \eqref{regou2}.

Let us consider, instead, the case when \eqref{gl2} holds.
Under this condition, using \eqref{regou}, \eqref{regochi} and
Remark~\ref{suregog}, we modify
\eqref{stimaF} as
%
\beeq{stimaF2}
 0\le\frac14\|u\|_{\calV}^2
  +\frac12|\chi_t|^2
  \le \calZ:=\calF+\|g\|_{\calV^*}^2
  \le \|u\|_{\calV}^2
  +\frac12|\chi_t|^2
  +\frac32\|g\|_{\calV^*}^2.
\end{equation}
Thus, in view of the first line of \eqref{conto45}, it is not difficult
to deduce
\beeq{conto45.2}
  \ddt \calZ
   \le c\big(\calZ^2+1+\|g\|_{\calV^*}^2+\|g_t\|_{\calV^*}^2\big)
    \quext{a.e.~in }\,(0,\infty),
\end{equation}
whence the convergence $u(t)\to 0$ with respect to
the {\sl strong}\/ topology on $\calV$ is a consequence of
\cite[Lemma~6.2.1, p.~225]{Zh}. Then, recalling \eqref{nortetau},
in this case it happens that \eqref{utetainfi} is improved into
$$
   u(t)\to 0,~~\teta(t)\to\teta_\infty
    \quext{strongly in }\,\calV.
$$

To conclude, it remains to show that the limit value $\chi_\infty$ in
\eqref{chiinfi} satisfies \eqref{prostaz}, and this
can be done in a completely standard way.
Namely, defining $(\chi_n,u_n):(0,1)\to V\times\calV$ as
$(\chi_n,u_n)(\cpt):=(\chi,u)(t_n+\cpt)$, from \eqref{regochi} and
\eqref{regou} it is clear that
\beeq{new-2}
(\chi_n,u_n)\to (\chi_\infty,0) \quext{strongly in }\, C^0([0,1]; H) \times L^2(0,1;\calV) .
\end{equation}
as $n \nearrow \infty$. Here, of course we also used the strong convergence
of $\chi_{n,t}$ to $0$ in $L^2(0,1;H)$.
Then, by passing to the limit in
\beeq{phasen}
    \chi_{n,t} + A \chi_n + W'(\chi_n)= \lambda'(\chi_n )u_n
\quext{in }\,V^*, \ \hbox{ a.e. in } \, (0,1),
\end{equation}
it is not difficult to check that $\chi_\infty$ solves the
stationary problem. Indeed, owing to \eqref{regochi2}, weak star compactness,
and \eqref{new-2} it turns out that $\chi_n \to \chi_\infty$ weakly star in
$L^\infty(0,1;H^2(\Omega))$, whence $A\chi_n \to A\chi_\infty$ weakly star in
$L^\infty(0,1;H)$. Moreover, using \eqref{regoW2} and exploiting
the maximal monotonicity of $W'+\kappa\Id$ (which is a continuous and
increasing function
thanks to \eqref{hpW}), one verifies that $W'(\chi_n)$ tends to $W'(\chi_\infty)$
weakly star in $L^\infty(0,1;H)$ (and then weakly in $L^2(\Omega \times (0,1))$)
with the help of the strong convergence $\chi_n  \to \chi_\infty$ in $L^2(\Omega
\times (0,1))$ and of \cite[Prop.~1.1, p.~42]{Ba}. Finally, thanks to \eqref{new-2},
\eqref{hplambda} and the continuous embedding $V\subset L^6(\Omega)$,
we infer that
$$
  \lambda'(\chi_n)u_n\to 0
   \quext{strongly in }\,L^2(0,1;L^{3/2}(\Omega)),
$$
and consequently the right hand side of \eqref{phasen} tends to $0$
in $L^2(0,1;V^*)$. This concludes the proof of Theorem~\ref{teoomega1}.
\dimbox

\vspace{2mm}

\noindent%
{\bf Proof of Theorem~\ref{teoomega2}.}~~%
We proceed partly as in \cite[Sec.~3]{CJ}, \cite{GPS2}.
Let us assume, for simplicity, that \eqref{condelta1} holds
(otherwise, we can replace $\zeta$ with
a value $\zeta_0$ such
that \eqref{condelta2} is satisfied and notice that
Theorem~\ref{teoloj} still holds with $\zeta_0$ in
place of $\zeta$ in \eqref{disloj}).
Letting $\chi_\infty$ be an element of the $\omega$-limit
of $\chi(\cdot)$, we can set (cf.~Theorem~\ref{teoloj} for the
notation)
\beeq{jen1}
  \Sigma:=\big\{
    t>0:\|\chi(t)-\chii\|_{V\cap C^0(\barO)}\le \epsilon/3\big\}.
\end{equation}
Clearly, $\Sigma$ is unbounded. Next,
for $t\in\Sigma$, we put
\beeq{jen2}
  \tau(t):=\sup\big\{
    t'\ge t:\sup_{s\in[t,t']}\|\chi(s)-\chii\|_{V\cap C^0(\barO)}
     \le \epsilon\big\}
\end{equation}
and observe that, by continuity, $\tau(t)>t$ for all
$t\in\Sigma$. Let us fix $t_0\in\Sigma$ and
divide ${\cal J}:=[t_0,\tau(t_0))$ (where
$\tau(t_0)$ might well be $+\infty$)
into two subsets:
\beal{A1}
  & A_1:=\Big\{
   t\in\calJ:|\chi_t(t)|+\|u(t)\|_{\calV}>\Big(\int_t^{\tau(t_0)}
   \|g(s)\|_{\calV^*}^2\,\dis\Big)^{1-\zeta}\Big\},\\
 \label{A2}
  & A_2:=\calJ\setminus A_1.
\end{align}
Next, 
we define (cf.~\eqref{enloj})
\beeq{defPhi}
  \Phi(t):=\io j(\teta(x,t))\,\dix
   +\frac12\int_t^{\tau(t_0)} \|g(s)\|_{\calV^*}^2\,\dis
   +E(\chi(t))-E(\chi_\infty).
\end{equation}
Then, it is not difficult to see that
\beeq{jen3}
  \Phi'(t)\le
   -\Big(
     \frac12\|u(t)\|_{\calV}^2
     +|\chi_t(t)|^2\Big).
\end{equation}
We remark that $\Phi$ is absolutely continuous
thanks to \eqref{regoteta2}--\eqref{regoW2} and
\cite[Lemme~3.3, p.~73]{Br}. This justifies the
above computation for a.e.~$t\in\calJ$.
We then have (cf.~also \cite[(3.2)]{Je})
\beeq{jen4}
  \ddt\big(|\Phi|^{\zeta}\sign\Phi\big)(t)\le
   -\zeta|\Phi(t)|^{\zeta-1}
    \Big(
     \frac12\|u(t)\|_{\calV}^2
     +|\chi_t(t)|^2\Big).
\end{equation}
Now, let us estimate $\Phi$ from above.
If $t\in A_1$, thanks to \eqref{enloj}
and Theorem~\ref{teoloj}, we obtain
\bealo
  |\Phi(t)|^{1-\zeta}
   & \le|E(\chi(t))-E(\chii)|^{1-\zeta}
   +\Big|\io j(\teta(t))\Big|^{1-\zeta}
    +\Big|\int_t^{\tau(t_0)}\|g(s)\|_{\calV^*}^2\,\dis\Big|^{1-\zeta}\\
 \label{conto71}
  & \le c_\ell\|-\chi_t(t)+\lambda'(\chi(t))u\|_{V^*}
   +\Big|\io j(\teta(t))\Big|^{1-\zeta}
   +|\chi_t(t)|+\|u(t)\|_\calV,
\end{align}
where we also used \eqref{A1}.
Note now that, by \eqref{hplambda}, the
last of \eqref{regochi2}, and well-known
continuous embeddings, we have
\beeq{conto72}
  c_\ell\|-\chi_t(t)+\lambda'(\chi(t))u\|_{V^*}
    \le c\big(|\chi_t(t)|+\|u(t)\|_{\calV}\big).
\end{equation}
Moreover, by convexity of $j$, $j(\teta_\infty)=0$,
and H\"older's inequality, we infer
\beeq{conto73}
  0 \le\io j(\teta(x,t))\,\dix
  \le\io j'(\teta(x,t))(\teta(x,t)-\teta_\infty)\,\dix
   \le |u(t)||\teta(t)-\teta_\infty|.
\end{equation}
Using once more \eqref{nortetau} together with
\eqref{regou2} and recalling that $\zeta\in(0,1/2)$,
one gets
\beeq{conto74}
  \Big|\io j(\teta(t))\Big|^{1-\zeta}
  \le c\|u(t)\|_{\calV}^{2(1-\zeta)}
  \le c\|u(t)\|_{\calV}.
\end{equation}
Thus, collecting \eqref{conto71}--\eqref{conto74}, we finally
have
\beeq{conto75}
  |\Phi(t)|^{1-\zeta}
   \le c\big(|\chi_t(t)|+\|u(t)\|_\calV\big),
\end{equation}
whence from \eqref{jen4} we obtain
\beeq{jen5}
   \|u(t)\|_{\calV}+|\chi_t(t)|
  \le -\frac c {4\zeta}\ddt\big(|\Phi|^{\zeta}\sign\Phi\big)(t).
\end{equation}
Since $\Phi$ is decreasing by \eqref{jen3},
integration in time entails that
$\|u\|_{\calV}$ and $|\chi_t|$ are summable over $A_1$.
Of course, the same holds over $A_2$ by \eqref{effe2}
and \eqref{A1}--\eqref{A2}.
Thus, we conclude that $\chi_t\in L^1(\calJ;H)$.

From this point on, the proof proceeds exactly
as in \cite[Sec.~3]{GPS2}. Namely, a simple
contradiction argument
yields that $\tau(t_0)=\infty$ as
$t_0\in\Sigma$ is sufficiently large. This implies that
$\chi_t\in L^1(t_0,+\infty;H)$, whence the convergence
(in $H$) of the whole trajectory
$\chi(t)$ to $\chii$ follows. More
precisely, this convergence holds strongly in
$V\cap C(\overline\Omega)$ by precompactness of the
trajectory (cf.~\eqref{regochi2}).
Finally, the technical
argument of \cite[Sec.~3]{GPS2} leading to estimate \eqref{rate1}
(or \eqref{rate2}) can be repeated just by adapting
the notation.\dimbox
%
%



\end{document}